\documentclass{article}
\usepackage[utf8]{inputenc}

\usepackage{amsmath}
\usepackage{amssymb}
\usepackage{amsthm}
\usepackage{amsfonts}
\usepackage{graphicx}
\usepackage{multirow}
\usepackage{xcolor}
\usepackage{url}

\usepackage{subcaption}

\newtheorem{theorem}{Theorem}
\newtheorem{lemma}{Lemma}
\newtheorem{assumption}{Assumption}
\newtheorem{corollary}{Corollary}
\newtheorem{remark}{Remark}

\newcommand{\labeltheoone}{1}
\newcommand{\labeltheotwo}{2}
\newcommand{\labeltheothree}{3}
\newcommand{\labeltheofour}{4}
\newcommand{\labellemone}{1}

\usepackage{natbib}

\usepackage{xcolor}

\usepackage{authblk}

\title{Dimension estimation in PCA model using high-dimensional data augmentation}

\author{U. Radoji\v{c}i\'c}
\affil{Institute of Statistics and Mathematical Methods in Economics, Vienna University of Technology, Austria}
\author{J. Virta}
\affil{Department of Matematics and Statistics, University of Turku, Finland}

\date{}

\begin{document}

\maketitle

\begin{abstract}
	We propose a modified, high-dimensional version of a recent dimension estimation procedure that determines the dimension via the introduction of augmented noise variables into the data. Our asymptotic results show that the proposal is consistent in wide high-dimensional scenarios, and further shed light on why the original method breaks down when the dimension of either the data or the augmentation becomes too large. Simulations are used to demonstrate the superiority of the proposal to competitors both under and outside of the theoretical model.
\end{abstract}

\section{Introduction}\label{sec:intro}

In this work, we revisit the classical problem of estimating the latent dimension in principal component analysis. Numerous solutions to this problem have been proposed in the literature, see, e.g., \cite{luo2016combining, nordhausen2021asymptotic, bernard2024power, virta2024robust} for some recent works. The standard solutions are predominantly based on sequential subsphericity testing, information-theoretic criteria, or risk minimization.


A different approach to dimension estimation is taken in a recent proposal known as predictor augmentation \citep{luo2021order}. The full description of the method is given in Section~\ref{sec:setting}, but on a heuristic level, in predictor augmentation the observed $n \times p$ data is augmented into a sample of size $n \times (p + r)$, where $r$ is essentially a tuning parameter and the added $nr$ variables are drawn i.i.d. from a normal distribution with a specific variance. The purpose behind the augmentation is that the added variables, being pure noise, get mixed with the actual noise subspace, allowing one to pinpoint the jump from the signal subspace to the noise subspace in the spectrum of the covariance matrix of the augmented data.

The outcome of the predictor augmentation is a function $\phi_n: \{0, \ldots, p\} \to \mathbb{R}$, whose minimizer $d_n$ is taken as the estimate of the true latent dimension $d$. In \cite[Theorem~5]{luo2021order} it is shown that this estimator is, for every fixed $r$, consistent as $n \rightarrow \infty$ under certain mild technical conditions. Notably, this result was given under the assumption of finite $p$ and it turns out that the consistency can fail in high-dimensional scenarios where the dimension $p \equiv p_n$ and/or the number of augmentations $r \equiv r_n$ are allowed to grow with $n$. 

To illustrate this phenomenon, we conducted an empirical study detailed in the supplementary Appendix A. 
The results highlight that in the ``low-dimensional'' scenario ($p_n = 0.025n$, $r_n = 0.01n$), the correct estimate is obtained. However, the other cases lead to misestimation. For example, when $p_n = 0.25n$, the choices $r_n = 0.01n$ and $r_n = 0.5n$ lead to over- and underestimation, respectively, indicating that the proper choice of $r_n$ lies somewhere in between. Even when the dimensionality is small ($p_n = 0.025n$), taking $r_n$ too large can lead to failure, contrasting with intuition from \cite{radojicic2022}, where larger $r_n$ is generally deemed beneficial for fixed $p$.

The reason for the previous issues is that when either the data dimension $p_n$ or the augmentation dimension $r_n$ is comparable in magnitude to $n$, we enter the domain of high-dimensional asymptotics and the classical arguments in \cite{luo2021order, radojicic2022} stop working. As such, the purpose of the current work is three-fold: (i) We investigate what causes the inconsistency of the augmentation estimator by a careful study of its asymptotic properties in the doubly high-dimensional setting where both the data dimension $p_n$ and the augmentation dimension $r_n$ diverge to infinity in the sense that $(p_n/n, r_n/n) \rightarrow (\gamma_p, \gamma_r) > 0$. As one consequence of our results, we show that for every $\gamma_p$, there exist augmentation rates $\gamma_r$, which lead to inconsistent predictor augmentation. (ii) We use our asymptotic results to derive a corrected estimator that is consistent in high-dimensional settings under very mild conditions on the signal strength. (iii) We use simulations to compare our estimator to both the original predictor augmentation and the high-dimensional subsphericity-based estimator by \cite{schott2006high}. The results indicate that our proposal surpasses both competitors and tolerates deviations from the model exceedingly well.


The paper is organized as follows. In Section \ref{sec:setting}, we review the original predictor augmentation estimator by \cite{luo2021order} and introduce our mathematical framework. 
\textcolor{black}{In Section \ref{sec:asymp_known_noise}, we derive the high-dimensional asymptotic behavior of the predictor augmentation components, which are then used in Section \ref{sec:inconsistency} to identify scenarios of inconsistency.} 
We propose our modified estimator and establish its consistency in Section \ref{sec:our_proposal}. Simulation results are presented in Section~\ref{sec:simulations}. Additional numerical experiments and the proofs of the technical results are given in the supplementary Appendices A and B, respectively.


\section{Model and the estimator}\label{sec:setting}

We assume throughout that $x_{1, n}, \ldots, x_{n, n}$ is a random $n$-indexed sample (triangular array) from the $p_n$-variate normal distribution $\mathcal{N}_{p_n}(\mu_n, \Sigma_n)$ where the $p_n$ eigenvalues of the covariance matrix $\Sigma_n$ are assumed to be $\lambda_1 + \sigma^2 > \cdots > \lambda_d + \sigma^2 > \sigma^2 = \cdots = \sigma^2$. That is, the constants $\lambda_1, \ldots, \lambda_d, \sigma^2$ and the true signal dimension $d \geq 1$ do not vary with $n$, and increasing $n$ simply has the effect of adding more noise dimensions in the model. Such models are commonly known as spiked covariance structures; see \cite{yao2015large}. The assumption that the spike eigenvalues are distinct is made for mathematical convenience. Our main objective in this paper is to, given the sample $x_{1, n}, \ldots, x_{n, n}$, estimate the signal dimension $d$.

To estimate $d$ with predictor augmentation \citep{luo2021order}, the observations $x_{1, n}, \ldots, x_{n, n}$ are first augmented with random noise simulated from normal distribution. That is, we form the augmented sample $z_{i, n} = (x_{i, n}', \sigma_n s_{i, n}')'$, $i = 1, \ldots, n$, where $s_n \sim \mathcal{N}_{r_n}(0, I_{r_n})$, and $\sigma_n$ is an estimator of the noise standard deviation $\sigma$. The dimensionality $r_n$ of the augmentation is a user-specified tuning parameter. Denoting by $S_n$ the sample covariance matrix of the $(p_n + r_n)$-dimensional augmented sample, we then compute a set $u_{1, n}, \ldots, u_{p_n, n}$ of any of its leading $p_n$ eigenvectors. Decomposing the eigenvectors as $u_{j, n} = (u_{j, n, A}', u_{j, n, B}', u_{j, n, C}')'$ where the dimensionalities of the three parts are $d, p_n - d, r_n$, respectively, the predictor augmentation estimator of $d$ is then based on the sequence of squared norms $\| u_{1, n, C} \|^2, \ldots, \| u_{p, n, C} \|^2$. Heuristically, a jump is seen in the magnitudes of these norms when we cross from the signal subspace into the noise subspace. \cite{luo2021order} further combine the eigenvector information with a standardized scree plot computed from the $p_n + 1$ first eigenvalues $\tau_{1, n}, \ldots, \tau_{p_n + 1, n}$ of $S_n$, and define the function $\phi_n: \{0, \ldots, p_n \} \to \mathbb{R}$, acting as
\begin{align}\label{eq:luoliobjective}
    \phi_n(k) = \sum_{j = 0}^k \| u_{j, n, C} \|^2 + \frac{\tau_{k + 1, n}}{1 + \sum_{j = 1}^{k + 1} \tau_{j, n}},
\end{align}
where we take $\| u_{0, n, C} \|^2 = 0$. The estimator of the signal dimension $d$ is then obtained as the minimizer of $\phi_n$. Additionally, \cite{luo2021order} considered independently repeating the augmentation several times and averaging the eigenvector norms over these. We ignore this possibility here, as it does not affect the limiting values of the studied quantities.

In the following sections, we study the behavior of the above estimator under the assumption that $p_n/n \rightarrow \gamma_p \in (0, \infty)$ and $r_n/n \rightarrow \gamma_r \in (0, \infty)$, as $n \rightarrow \infty$, for some constants $\gamma_p, \gamma_r$. Throughout the work, we make the following assumption.

\begin{assumption}\label{assu:main_assumption}
    The smallest signal eigenvalue satisfies $\lambda_d > \sigma^2 \sqrt{\gamma_p + \gamma_r}$.
\end{assumption}

Assumption \ref{assu:main_assumption} is connected to the so-called Baik-Ben Arous-P\'ech\'e transition which states that the sample estimate of any signal eigenvalue in the interval $(\sigma^2, \sigma^2(1 + \sqrt{\gamma_p + \gamma_r}))$ asymptotically behaves as it was a noise eigenvalue, the behavior extending also to the corresponding eigenvector \citep{yao2015large}. Hence, Assumption \ref{assu:main_assumption} is natural and essentially necessary, as without it the effective dimension of our model would actually be smaller than $d$. Similar assumptions are commonly made, see, e.g., \cite{alaoui2020fundamental, jagannath2020statistical, mukherjee2023consistent}.

\section{Asymptotics of predictor augmentation}\label{sec:asymp_known_noise}

For mathematical convenience, we make the simplifying assumption that the noise variance $\sigma^2$ is known. This allows us to bypass its estimation and use the true value of $\sigma$ in the generation of the augmented variables. The simulations in Section \ref{sec:simulations} later show that the practical behavior of the method changes remarkably little when $\sigma$ is replaced with a consistent estimator $\sigma_n$ of it.

We begin by establishing the probability limits of the squared norms $\| u_{j, n, C} \|^2$ of the augmented parts of the eigenvectors.

\begin{theorem}\label{theo:known_noise_1}
    Fix a constant $K \in \mathbb{N}$ such that $K > d$. Under Assumption \ref{assu:main_assumption}, we have, as $n \rightarrow \infty$,
    \begin{align*}
        \| u_{j, n, C} \|^2 &\rightarrow_p \frac{\gamma_r \sigma^2}{\lambda_j} \left( \frac{\lambda_j + \sigma^2}{\lambda_j + (\gamma_p + \gamma_r) \sigma^2} \right), \quad j = 1, \ldots, d, \\
        \| u_{j, n, C} \|^2 &\rightarrow_p \frac{\gamma_r}{\gamma_p + \gamma_r}, \quad j = d + 1, \ldots, K.
    \end{align*}
\end{theorem}

\begin{remark}\label{rem:generalization}
    Theorem \ref{theo:known_noise_1} could still be generalized to have $\min \{ p_n + r_n, n - 1 \}$ in place of the constant $K$, assuming that $\gamma_p + \gamma_r \neq 1$, see \cite[Theorem 2.17]{Knowles2014} for details.
\end{remark}

To accompany the norms, the augmentation estimator also requires estimating the signal eigenvalues. \cite{luo2021order} estimate them as the largest $p_n +~1$ eigenvalues $\tau_{j, n}$, $j = 1, \ldots, p_n +~1$ of the sample covariance matrix of the augmented sample $z_{i, n}$. The following result gives the probability limits of the eigenvalues $\tau_{j,n}$ illustrating further the inconsistency of the traditional estimators in the high-dimensional regime. Here  $F_\gamma^{-1}$ denotes the quantile function of the Marchenko-Pastur distribution with the concentration parameter $\gamma > 0$.

\begin{theorem}\label{theo:known_noise_2}
    Under Assumption \ref{assu:main_assumption}, we have for the signal eigenvalues that
    \begin{align*}
       \tau_{j, n} &\rightarrow_p (\lambda_j + \sigma^2) \{ 1 + (\gamma_p + \gamma_r) \sigma^2 / \lambda_j \}, \quad j = 1, \ldots, d,
    \end{align*}
    Whereas, for the noise eigenvalues we have the following two cases,
    \begin{enumerate}
        \item[(i)] If $\gamma_p + \gamma_r \in (0, 1]$, then for a sequence $j_n > d$, such that $j_n/(p_n + r_n) \rightarrow 1 - q \in [0, 1]$ as $n \rightarrow \infty$, we have $\tau_{j_n, n} \rightarrow_p \sigma^2 F^{-1}_{\gamma_p + \gamma_r}(q)$.
        \item[(ii)] If $\gamma_p + \gamma_r \in (1, \infty)$, then for a sequence $j_n \in [d + 1, n - 1]$, such that $j_n/(n - 1) \rightarrow 1 - q \in [0, 1]$ as $n \rightarrow \infty$, we have $\tau_{j_n, n} \rightarrow_p \sigma^2 (\gamma_p + \gamma_r) F^{-1}_{(\gamma_p + \gamma_r)^{-1}}(q)$, whereas, $\tau_{j, n} = 0$ for all $j \geq n$.  
    \end{enumerate}
    
\end{theorem}

\begin{remark}\label{rem:bbp_effect}
    Assume that, contrary to Assumption \ref{assu:main_assumption}, we would have a very weak signal, in the sense that $\lambda_d = \sigma^2 \sqrt{\gamma_p + \gamma_r}$. In this case, arguing as in the proofs of Theorems \ref{theo:known_noise_1} and \ref{theo:known_noise_2}, one could show that both $\| u_{d, n, C} \|^2$ and $\| u_{d + 1, n, C} \|^2$ would then converge to the same constant, as would $\tau_{d, n}$ and $\tau_{d + 1, n}$. Hence, this further demonstrates our earlier discussion that, without Assumption~\ref{assu:main_assumption}, no method based on the asymptotic limits of the eigenvalues and eigenvectors of the sample covariance matrix is able to estimate $d$. 
\end{remark}

In the sequel, we will denote the probability limits of the eigenvector norms $\| u_{j, n, C} \|^2$ and the eigenvalues $\tau_{j, n}$ by $\| u_{j, C} \|^2$ and $\tau_{j}$, respectively. Moreover, the corresponding point-wise limit of the objective function $\phi_n$ in \eqref{eq:luoliobjective} will be denoted by $\phi$, i.e., $\phi_n(k)\to_p \phi(k)$ for $k\geq 0$.

\section{When is predictor augmentation inconsistent?}\label{sec:inconsistency}

Essentially, Theorems \ref{theo:known_noise_1} and \ref{theo:known_noise_2} can be used to determine whether any given high-dimensional scenario leads to a consistent augmentation estimate of $d$. For example, plugging in the specification of the simulation scenario in Section \ref{sec:intro} and Appendix A with $d = 1$ reveals that, indeed, $\phi(1) > \phi(2)$ holds when $\gamma_p = 0.25$ and $\gamma_r = 0.01$, leading to the observed inconsistent estimate. Due to the complicated form of the objective function \eqref{eq:luoliobjective}, the full answer to the question ``for which parameter values is predictor augmentation consistent'' is not feasible to obtain, but several qualitative statements can be constructed, as formalized in the following theorem.
\begin{theorem}\label{thm:LuoLi_inconsistency}
    Under Assumption \ref{assu:main_assumption}, the following statements hold.
    \begin{itemize}
        \item[(i)] Given any $\lambda_d > 0$ and $\gamma_p > 0$, there exists $\gamma_r^0 \in (0, \lambda_d^2\sigma^{-4} - \gamma_p )$, such that $\phi(d) > \phi(d+1)$ for every $\gamma_r \in (0, \gamma_r^0)$.
        \item[(ii)] Given any $\gamma_p > 0$, there exists $\lambda_d > \sigma^2\sqrt{\gamma_p}$ small enough, such that for all $\gamma_r \in (0, \lambda_d^2\sigma^{-4} - \gamma_p )$ we have $\phi(d) > \phi(d+1) $.
        \item[(iii)]  Given any $\lambda_d > 0$, there exists $\gamma_p \in (0, \lambda_d^2\sigma^{-4})$ large enough, that for all $\gamma_r \in (0, \lambda_d^2\sigma^{-4} - \gamma_p )$ we have $\phi(d) > \phi(d+1) $.
    \end{itemize}
\end{theorem}


Statement (i) of Theorem \ref{thm:LuoLi_inconsistency} shows that there always exists a sufficiently small $\gamma_r$ for which the augmentation estimator is inconsistent, implying that $r_n$ must be large enough relative to $p_n$; in high-dimensional settings, using a finite number of augmentations is always insufficient for consistency. Statements (ii) and (iii) describe more challenging cases: (ii) asserts that for any data collection rate $\gamma_p$, there exists a scenario where the signal is asymptotically distinguishable from noise, yet the augmentation estimator remains inconsistent. Similarly, (iii) states that no matter how strong the signal-to-noise ratio is, there will always exist a sufficiently large $\gamma_p$ for which the signal remains distinguishable, but the augmentation still yields an inconsistent estimate.

\begin{figure}[t]
    \centering
\includegraphics[width=0.8\linewidth]{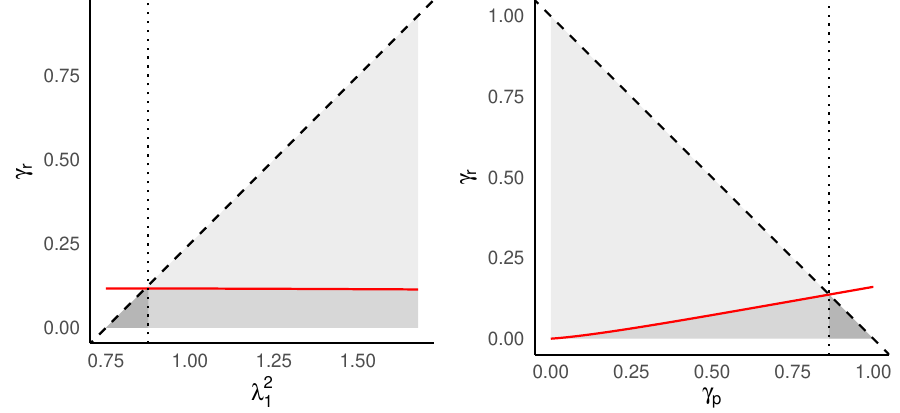}
    \caption{Graphical illustration of the results of Theorem \ref{thm:LuoLi_inconsistency}.}
    \label{fig:inconsitency_region}
\end{figure}

Figure \ref{fig:inconsitency_region} illustrates these inconsistency regions. In the left plot, for $\gamma_p=0.75$, $\sigma^2=1$, and $d=1$, the light gray region bounded by the dashed line represents $\gamma_r + \gamma_p < \lambda_1^2$, defining feasible values of $(\lambda_1, \gamma_r)$ given $\gamma_p=0.75$ (see the interpretation after Assumption \ref{assu:main_assumption}). The red curve shows, for each $\lambda_1$, the value of $\gamma_r^0$ from Theorem \ref{thm:LuoLi_inconsistency} (i) where $\phi(1) = \phi(2)$. The darker gray area below this curve indicates that even with a sufficiently strong signal (in the sense of Assumption~\ref{assu:main_assumption}), inconsistency can happen for every feasible $\gamma_r$ (i.e., $\gamma_r$ satisfying $\lambda_1^2 > \gamma_p + \gamma_r$). The right plot shows similar inconsistency regions for $\lambda_1=1$, $\sigma^2=1$, and $d=1$.

\section{High-dimensionally consistent estimator of $d$}\label{sec:our_proposal}

We next define an alternative to the augmentation function $\phi_n$ that retains its consistency for high-dimensional data. Our proposal has three key elements: (i) Unlike \cite{luo2021order}, we estimate the signal eigenvalues directly from the original data (instead of from the augmented data). This leads to less biased estimates of $\lambda_j$. (ii) We further correct for the remaining bias in the estimated eigenvalues by applying an appropriately chosen transformation to them. (iii) We combine the eigenvalue and eigenvector information not by summing them, but by adjusting the norms $\| u_{j, n, C} \|^2$ such that the jump from the signal to the noise becomes apparent in their plot.

As per item (i) above, we assume throughout this section that the eigenvalues $\tau_{j, n}$ have been estimated from the sample covariance matrix of the original data. In this case the equivalent of Theorem \ref{theo:known_noise_2} holds with $\gamma_r = 0$. Define next the debiasing function $f_n:\mathbb{R} \to \mathbb{R}$ as
\begin{align*}
    f_n(\tau) = \frac{1}{2} \{ \tau - \sigma^2 (1 + p_n/n) \} + \frac{1}{2} [ \{ \tau - \sigma^2 (1 + p_n/n) \}^2 - 4 \sigma^4 (p_n/n) ]_{+}^{1/2},
\end{align*}
where $[ a ]_{+} := \max\{0 ,a\}$ is the soft-thresholding function. By Theorem~\ref{theo:known_noise_2} and direct computation, we see that $f_n(\tau_{j, n}) \rightarrow_p \lambda_j$ for all $j = 1, \ldots, d$, showing that $f_n$ allows for the unbiased estimation of the spikes. Using the debiased estimates $f_n(\tau_{j, n})$ and the eigenvector norms, we then construct the following quantities,
\begin{align}\label{eq:hjn}
    h_{j, n} := \frac{f_n(\tau_{j, n}) \{ f_n(\tau_{j, n}) + (p_n/n + r_n/n) \sigma^2 \}}{\{ f_n(\tau_{j, n}) + \sigma^2 \} } \| u_{j, n, C} \|^2.
\end{align}
Fix now a large constant $K \in \mathbb{N}$. We propose estimating the signal dimension as that value $j \in 1, \ldots, K$, for which the successive difference $h_{j + 1, n} - h_{j, n}$ is minimized. Our next result implies that this estimator is consistent under all possible high-dimensional regimes $\gamma_p, \gamma_r > 0$, as long as the signal is identifiable in the sense of Assumption \ref{assu:main_assumption}. 

\begin{theorem}\label{theo_our_consistency}
        Fix any $K \in \mathbb{N}$ such that $K > d$. Then, under Assumption \ref{assu:main_assumption}, we have that
     \begin{align*}
         h_{j + 1, n} - h_{j, n} \rightarrow_p 0,
     \end{align*}
     for all $j = 1, \ldots, K$ such that $j \neq d$. Moreover,
     \begin{align*}
         h_{d + 1, n} - h_{d, n} \rightarrow_p \frac{- \sigma^2\gamma_r^2}{(\sqrt{\gamma_p} + 1)(\gamma_p + \gamma_r)} < 0.
     \end{align*}
\end{theorem}

\begin{corollary}\label{cor:cor1}
    Fix any $K \in \mathbb{N}$ such that $K > d$. Then, under Assumption \ref{assu:main_assumption}, we have that $d_n := \textrm{argmin}_{j = 1, \ldots , K} \{ h_{j + 1, n} - h_{j, n} \}$ satisfies $d_n \rightarrow_p d$. 
\end{corollary}


\begin{remark}
    In Theorem \ref{theo_our_consistency} we search the true signal dimension within the finite set $\{1, \ldots, K\}$. This is practical since, in dimension reduction, one typically assumes that the signal of the data is captured by a relatively small amount of factors. However, the search interval could also be widened and allowed to grow with $n$, in the same sense as discussed in Remark \ref{rem:generalization}.
\end{remark}

Theorem \ref{theo_our_consistency} quantifies the limiting "jump" in the successive differences of the adjusted augmented norms $h_{j,n}$, $j=1,\dots,K>d$, once the true latent dimension is reached. We can therefore use this result to "tune" the parameter $\gamma_r$, so that the jump is maximized. If $\gamma_p>0$, substituting $\gamma_r=c\gamma_p$ in the limit of $h_{d+1,n}-h_{d,n}$, we get
\begin{align}\label{eq:temp1}
     h_{d + 1, n} - h_{d, n} \rightarrow_p T(\gamma_p,c,\sigma) := \frac{- \sigma^2 c^2 \gamma_p }{(\sqrt{\gamma_p}+1)(1+c)} < 0.
\end{align}
As $c,\gamma_p>0$, it is easily shown that for every $\gamma_p,\sigma^2>0$, $c\mapsto T(\gamma_p,c,\sigma)$ is strictly decreasing in $c$, implying that $\gamma_r$ should be chosen as large as possible, taking care not to violate Assumption~\ref{assu:main_assumption}. If $\gamma_p=0$, then $|h_{d + 1, n} - h_{d, n}| \rightarrow_p \sigma^2\gamma_r$, which is again maximized by taking $\gamma_r$ as large as possible. Also, this illustrates the potential advantage of the proposed method with respect to the original predictor augmentation even in low-dimensional settings, where one would benefit from the ability to take $r_n$ as large as possible; see also the discussion in Section \ref{sec:intro}.

As in practice the noise variance $\sigma^2$ is unknown, We conclude the section with a result giving a number of possible consistent estimators of $\sigma^2$. 
\begin{lemma}\label{lemma:sigma_estimation}
     Let $\tau_{j_n,n}$ be the $j_n$th eigenvalue of the sample covariance of $x_{1,n},\dots,x_{n,n}$. Then for $j_n\in\{d+1,\dots,p_n\}$ such that $j_n/p_n\to 1-q\in [0,1]$, as $n\to \infty$, and under Assumption \ref{assu:main_assumption},
     $$
    \hat{\sigma}^2_{j_n} := g(\tau_{j_n,n})=\left.\begin{cases}
        \tau_{j_n,n}/F_{\gamma_p}^{-1}(q),& \gamma_p\leq 1,\\
        \tau_{j_n,n}/(\gamma_p F_{\gamma_p^{-1}}^{-1}(q)),&\gamma_p > 1,
    \end{cases}\right\}\quad\to_p\quad \sigma^2,
    $$
    where $F^{-1}_\gamma$ is as in Theorem~\ref{theo:known_noise_2}.
\end{lemma}
In simulations, we use the corrected median eigenvalue as the estimator, with $j_n = \lfloor p_n/2 \rfloor$; for large $n$, the assumption $j_n > d$ is trivially satisfied, as the signal dimension $d$ is fixed and finite.
\section{Simulations}\label{sec:simulations}

The goal of this simulation study is threefold. First, we evaluate the validity of the proposed high-dimensional predictor augmentation (HDPA) based on the consistency result in Corollary \ref{cor:cor1}. Second, we compare HDPA with two alternatives: the original predictor augmentation by \cite{luo2021order} (PA) and the high-dimensional subsphericity-based estimator by \cite{schott2006high}. Finally, we explore HDPA's sensitivity to violations of the Gaussianity assumption. The \texttt{R} code for reproducing the results presented in this section is available at \url{https://github.com/uradojic/High-dimensional-data-augmentation}. Data is generated according to two models:
\begin{align*}
    \text{Model 1:}\quad X\sim\mathcal{N}_{p_n}({0}_{p_n},{\Sigma}_{p_n})\,,\,\quad\quad\quad \text{Model 2:}\quad X={\Sigma}^{1/2}_{p_n}(2{B}-{1}_{p_n}),
\end{align*}
where ${B}=(B_1,\dots,B_{p_n})'$ has i.i.d. entries $B_i\sim\mathrm{Bernoulli}(0.5)$ and ${1}_{p_n}\in\mathbb{R}^{p_n}$ is a vector of ones. In both models the covariance matrix $\boldsymbol{\Sigma}_{p_n}$ is taken to be of the form ${\Sigma}_{p_n}=\mathrm{diag}(5,4.8,\dots,3.2,3,0,\dots,0) + {I}_{p_n}$, giving latent dimension $d=11$ and noise variance $\sigma^2=1$. 

As both PA and HDPA require the estimate of the noise variance, to assess the effect of the estimate of the noise variance for both estimators, we use the oracle $\sigma^2=1$ as well as the eigenvalue-based estimator; for augmentation-based estimators 
we use the corrected median eigenvalue estimator $\hat{\sigma}^2_{\lfloor n/2 \rfloor}$ defined in Lemma \ref{lemma:sigma_estimation}.



\begin{figure}[t!]
    \centering
    \includegraphics[width=\linewidth]{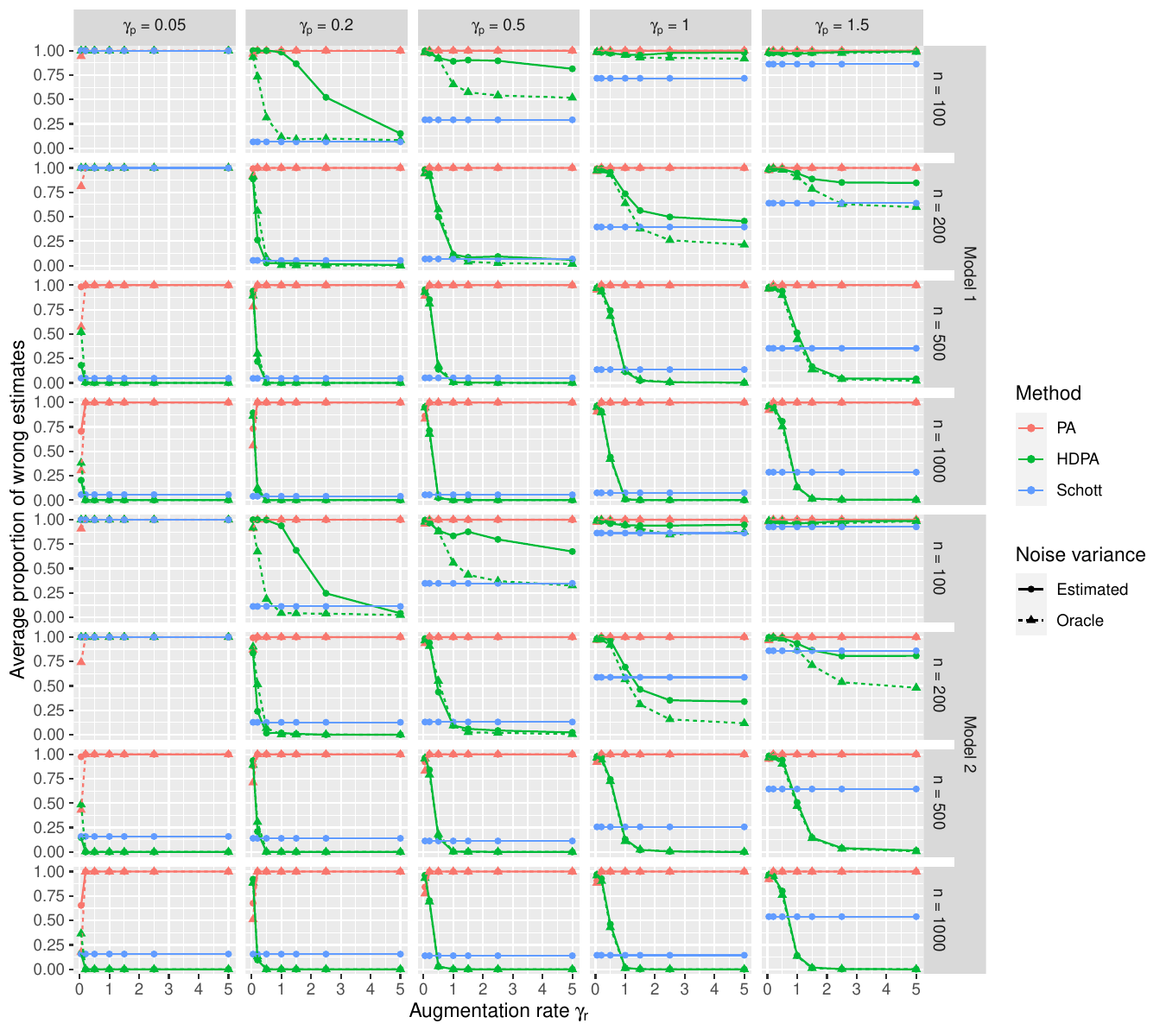}
    \caption{Average proportion of the wrong estimates in $1000$ replicates across $20$ different settings 
    for 
    Models 1 and 2.}
    \label{fig:M1}
\end{figure}

For each of two models, we replicate $m=1000$ data sets of size $n=100,\,200,\,500,\,1000$, and dimensionality $p_n=\gamma_pn$, for $\gamma_p=0.05,\,0.2,\,0.5,\,1,\,1.5$. To study the effect of the augmentation dimension $r_n$, we estimate the latent dimension in all settings using $r_n=\gamma_rn$, for $\gamma_r=0.05,\, 0.2,\, 0.5,\, 1,\, 1.5,\, 2.5,\, 5$. The average proportions of wrong estimates for each of the three estimators are given in Figure \ref{fig:M1} for both Models 1 and 2. 

The results show that HDPA significantly outperforms the original PA, which fails to estimate the latent dimension accurately and shows no improvement with larger sample sizes. This aligns with the introductory discussion and Theorem~\ref{thm:LuoLi_inconsistency}. For $n \geq 500$, HDPA achieves high accuracy across all settings and models. Given its reliance on the limiting behavior of sample covariance eigenvalues and eigenvectors, the need for a sufficiently large sample size is expected.

Additionally, Figure \ref{fig:M1} demonstrates HDPA's strong performance in Model~2, highlighting its robustness against deviations from the Gaussianity assumption. In contrast, the estimator by \cite{schott2006high}, which also assumes Gaussianity, fails for Bernoulli data and is consistently outperformed by HDPA across nearly all combinations of $n$ and $\gamma_p$. In Gaussian Model 1, Schott's method is limited by its reliance on successive hypothesis testing, which is constrained by the $0.05$ significance level we use. While it could theoretically be adjusted with a sample-dependent significance level, this would be impractical and offer no guarantees in finite-sample scenarios.

Furthermore, in accordance with Theorem \ref{theo_our_consistency} and the corresponding discussion, we observed that the performance of HDPA increases with the increase of $\gamma_r$, with the note that Assumption~\ref{assu:main_assumption} is satisfied for all the parameter combinations. Finally, for sample size large enough, we observe no difference in the performance of HDPA when the oracle (known) noise variance is replaced by its consistent estimator. In conclusion, for a large enough number of augmentations, HDPA outperforms the competitors in almost all considered settings.  

\subsection*{Acknowledgments}
The work of JV was supported by the Research Council of Finland (Grants 335077, 347501, 353769).
The work of UR was supported by the Austrian Science Fund (FWF), [10.55776/I5799].

\appendix

\section{Additional numerical experiments}

\begin{figure}[h!]
    \centering
    \includegraphics[width=0.9
    \linewidth]{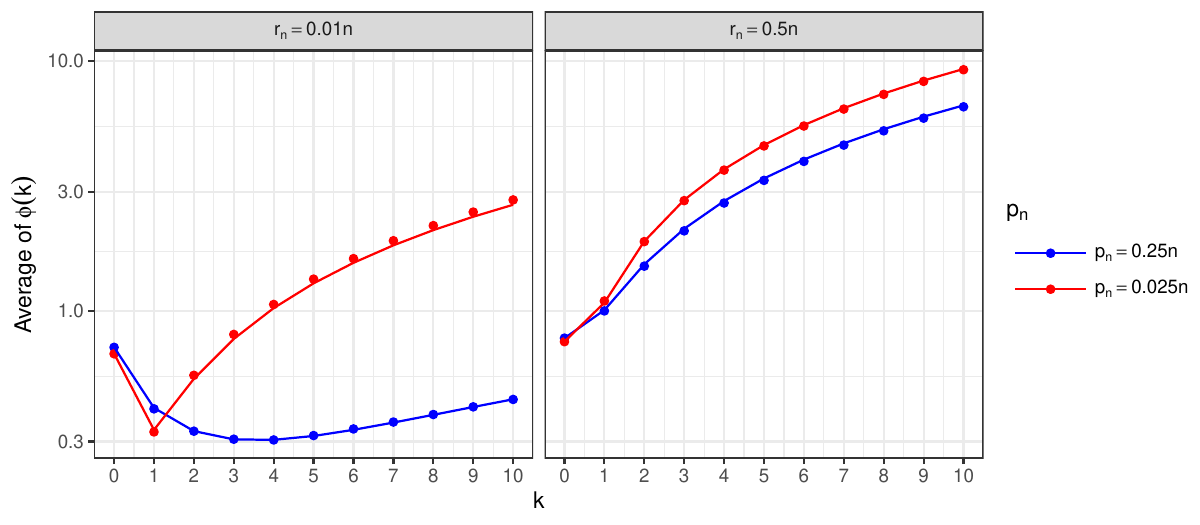}
    \caption{The average augmentation curves $\phi_n$ over $500$ replicates. 
    The solid lines indicate the 
    limiting values.}
    \label{fig:intro}
\end{figure}

To illustrate the inconsistency phenomena described in Section 1, we conducted an empirical study whose results are illustrated in Figure~\ref{fig:intro}. More specifically, Figure~\ref{fig:intro} shows the average augmentation curves $\phi_n$ (the points in the plot) over 500 independent replicates in four different settings with $n = 400$, $r_n = 0.01n, 0.5n$ and $p_n = 0.025n, 0.25n$. In all cases, the signal-to-noise ratio is $2$ and the true dimensionality is taken to be $d = 1$, meaning that optimally $\phi_n(k)$ should attain its minimum at $k = 1$. The solid lines in the plot show the corresponding pointwise limiting values of the objective function $\phi_n$ (when $n \rightarrow \infty$), computed based on the asymptotic results in our Section 3.

For visual convenience, the restrictions of the functions to the relevant part of the $x$-axis only are shown. From the plot, we observe that the correct estimate is indeed obtained in the ``low-dimensional'' case with $p_n = 0.025n, r_n = 0.01n$, but the three other cases lead to misestimation. E.g., when $p_n = 0.25$, the choices $r_n = 0.01n$ and $r_n = 0.5n$ lead to over and underestimation, respectively, indicating that the proper choice of $r_n$ is somewhere in between. We also observe, that even when the dimensionality of the data is small, $p_n = 0.025n$, the estimate can be ruined by taking $r_n$ too large. This goes against the intuition obtained in \cite{radojicic2022}, who generalize the approach from \cite{luo2021order} to $m$th order tensors, that larger $r_n$ should always be beneficial (for fixed $p$).

\section{Proofs of the main results}

\begin{proof}[Proof of Theorem \labeltheoone]
    We first note that since the sample covariance matrix $S_n$ involves centering, we may without loss of generality take $\mu = 0$. Then, denoting the eigendecomposition of the population covariance matrix by $\Sigma = O D O'$, we have $S_n = \mathrm{diag}(O D^{1/2}, \sigma I_{r_n}) W_n \mathrm{diag}(D^{1/2} O', \sigma I_{r_n})$ where $W_n$ is distributed as the sample covariance matrix of a sample of size $n$ from the $(p_n + r_n)$-variate standard normal distribution. By Theorem 3.4.4 (c) in \cite{mardia2024multivariate}, $W_n$ has the same distribution as $(n - 1)/n$ times the uncentered sample covariance matrix of a sample of size $n - 1$ from the $(p_n + r_n)$-variate standard normal distribution. This connection, together with the fact that $(n - 1)/n \rightarrow 1$, implies that we may also omit the sample centering from the calculation of $S_n$, and instead take the sample size to be $n - 1$ in the proofs.
    
    Next, since our main interest is on the joint distribution of $u_{1, n, C}, \ldots , u_{p, n, C}$, which is unaffected by transformations of the form $z_{i, n} \rightarrow \mathrm{diag}(V, \sigma I_{r_n}) z_{i, n}$ where $V$ is any $p_n \times p_n$ orthogonal matrix, it is without loss of generality that we may take $O = I_{p_n}$. Finally, we scale the whole data by $1/\sigma$, making the augmented data a sample (of size $n - 1$) from the multivariate normal distribution with the covariance matrix 
    \begin{align}\label{eq:sigma_n}
        \Sigma_n :=
        \begin{pmatrix}
        \Theta + I_d & 0 & 0 \\
        0 & I_{p_n - d} & 0 \\
        0 & 0 & I_{r_n}
        \end{pmatrix},
    \end{align}
    where $\Theta := \mathrm{diag}(\lambda_1/\sigma^2, \ldots, \lambda_d/\sigma^2)$.

    We first show the claim for $\| u_{j, n, C} \|^2$ where $j = 1, \ldots, d$, using the notation $m_{j, n} := (u_{j, n, B}', u_{j, n, C}')'$. Our setting is in the framework of Theorem 11.5 in \cite{yao2015large}, giving us
    \begin{align*}
        \| m_{j, n} \|^2 \rightarrow_P 1 - \alpha_j \frac{\psi'(\alpha_j)}{\psi(\alpha_j)},
    \end{align*}
    where $\psi, \psi'$ are as on page 202 of \cite{yao2015large} and $\alpha_j := 1 + \lambda_j/\sigma^2$. Plugging in their formulas, the probability limit of $\| m_{j, n} \|^2$ is seen to be
    \begin{align*}
        \| m_{j, n} \|^2 \rightarrow_P  \frac{(\gamma_p + \gamma_r) \sigma^2 (\lambda_j + \sigma^2)}{\lambda_j \{ \lambda_j + (\gamma_p + \gamma_r) \sigma^2 \} }.
    \end{align*}
    By the orthogonal invariance of the normal distribution, the form of the covariance matrix \eqref{eq:sigma_n} means that the distribution of $m_{j, n}$ is orthogonally invariant. Hence, $m_{j, n}/\| m_{j, n} \|$ is uniformly distributed in the unit sphere in $\mathbb{R}^{p - d + r_n}$ and $\| u_{j, n, C} \|^2 = \| m_{j, n} \|^2 Y_{j, n} $ where $Y_{j, n} \sim \mathrm{Beta}\{ r_n/2, (p - d)/2 \} \rightarrow_p \gamma_r/(\gamma_r + \gamma_p)$ as $n \rightarrow \infty$. Consequently,
\begin{align*}
    \| u_{j, n, C} \|^2 \rightarrow_p \frac{\gamma_r \sigma^2 (\lambda_j + \sigma^2)}{\lambda_j \{ \lambda_j + (\gamma_p + \gamma_r) \sigma^2 \} },
\end{align*}
for all $j = 1, \ldots, d$.

    For $j = d + 1, \ldots, p_n$, we have from \cite{Knowles2014}[Theorem 2.17] that $u_{j,n}'e_k \rightarrow_p 0$ for all $k = 1, \ldots, d$, which further implies that $\| u_{j, n, A} \| \rightarrow_p 0$, or conversely, $\| m_{j,n} \| \rightarrow_p 1$. Then, following the same argumentation as earlier, we obtain $\|u_{j,n,C}\|^2 \rightarrow_p \gamma_r/(\gamma_p + \gamma_r)$.
\end{proof}

\begin{proof}[Proof of Theorem \labeltheotwo]
    The result for the signal eigenvalues follows directly from Corollary 11.4 in \cite{yao2015large}. For the noise eigenvalues, we first observe that the eigenvalues of $\Sigma_n^{1/2} X'X \Sigma_n^{1/2}/(n - 1)$ are the same as the eigenvalues of $X'X \Sigma_n/(n - 1)$, where $\Sigma_n$ is as in \eqref{eq:sigma_n} and $X$ is a $(n - 1) \times (p_n + r_n)$ matrix full of standard normal variates. Theorem 2.14 \citep{yao2015large} now gives claim (i). For claim (ii), the final eigenvalues are equal to zero as the rank of $X'X \Sigma_n/(n - 1)$ is at most $n - 1$. For the non-trivial eigenvalues, we first observe that the non-zero spectrum of $X'X /(n - 1)$ is the same as that of $\{(p_n + r_n)/(n - 1)\} XX'/(p_n + r_n)$. Hence, Theorem 2.14 \citep{yao2015large} again gives the desired result, this time with inverse concentration.
\end{proof}

\begin{proof}[Proof of Theorem \labeltheothree]
   
    For $k\geq 0$, we let $\phi(k)$ to be the probability limit of $\phi_n(k)$, i.e. $\phi_n(k)\to_p \phi(k)$ for $k\geq 0$. Using the limiting results from Theorems \labeltheoone \ and \labeltheotwo , and denoting $\tau_i>0$ to be the probability limit of $\tau_{i,n}$, $i=1,\dots,d$, we obtain that for $\gamma_r,\gamma_p>0$, $\lambda_d \geq \sigma^2\sqrt{\gamma_p+\gamma_r}$,
    \begin{equation}\label{eq:phi_d-phi_(d-1)}
    \phi(d)-\phi(d+1)=-\frac{\gamma_r}{\gamma_r+\gamma_p}+\frac{\tau_{d+1}^2}{\left(1+\sum_{i=1}^{d+1}\tau_i\right)\left(1+\sum_{i=1}^{d+1}\tau_i+\tau_{d+1}\right)}.
    \end{equation}
    
    \textbf{Statement (i):} Using identity \eqref{eq:phi_d-phi_(d-1)} we conclude that 
    $$
 \phi(d)-\phi(d+1)>0 \iff -\gamma_r\left(1+\sum_{i=1}^{d+1}\tau_i\right)\left(1+\sum_{i=1}^{d+1}\tau_i+\tau_{d+1}\right)+(\gamma_r+\gamma_p)\tau_{d+1}^2>0.
    $$
    Function $g(\gamma_r)= -\gamma_r\left(1+\sum_{i=1}^{d+1}\tau_i\right)\left(1+\sum_{i=1}^{d+1}\tau_i+\tau_{d+1}\right)+(\gamma_r+\gamma_p)\tau_{d+1}^2$ is continuous in $\gamma_r$ and satisfies $g(0)>0$. Continuity of $g$ then implies that there exists $\gamma_r^0=\gamma_r^0(\gamma_p,\lambda_1,\dots,\lambda_d,\sigma^2)>0$ small enough such that 
    $g(\gamma_r)>0$ for every $\gamma_r {<} \gamma_r^0$, thus proving the statement (i).
    
    
    \textbf{Statement (ii):} Fix $\gamma_p>0$. Then, we will show that there exists $\varepsilon>0$ such that, for $\lambda_d=\sigma^2\sqrt{\gamma_p}+\varepsilon$, $\gamma_r<\varepsilon\sigma^{-4}(\varepsilon+2\sigma^2\sqrt{\gamma_p})$ implies that $\phi(d)-\phi(d+1)>0$. 
    For simplicity, we present the proof for $d=1$ and $\sigma^2=1$. The general case is proven analogously. Equation \eqref{eq:phi_d-phi_(d-1)} for $d=1$ gives
    \begin{equation}\label{eq:phi_1-phi_(2)}
     \phi(1)-\phi(2)>0 \iff \gamma_r\left(1+\tau_1+\tau_2\right)\left(1+\tau_1+2\tau_2\right)<(\gamma_r+\gamma_p)\tau_{2}^2.
\end{equation}
As $\tau_1\geq\tau_2\geq 1$ for any $\gamma_p,\gamma_r>0$ then 
\begin{equation*}
\gamma_r\left(1+3\tau_1\right)^2<(\gamma_r+\gamma_p)\Rightarrow \phi(1)-\phi(2)>0,    
\end{equation*}
i.e.,
\begin{equation*}
\gamma_r\left(\lambda_1+3(\lambda_1+1)(\lambda_1+\gamma_r+\gamma_p)\right)^2<\lambda_1^2(\gamma_r+\gamma_p)\Rightarrow \phi(1)-\phi(2)>0,  
\end{equation*}
since $\tau_1=(\lambda_1+1)(\lambda_1+\gamma_r+\gamma_p)\lambda_1^{-1}$. Substituting 
$\lambda_1=\sqrt{\gamma_p}+\varepsilon$ into the upper inequality we obtain
\begin{align}\label{eq:temp12}
\begin{split}
& \gamma_r \left(\sqrt{\gamma_p}+\varepsilon+3(\sqrt{\gamma_p}+\varepsilon+1 (\sqrt{\gamma_p}+\varepsilon+\gamma_r+\gamma_p)\right)^2<(\sqrt{\gamma_p}+\varepsilon)^2(\gamma_r+\gamma_p) \\
& \Rightarrow \phi(1)-\phi(2)>0.
\end{split}
\end{align}

For a fixed $\varepsilon > 0$, let now $0<\gamma_r<\varepsilon(\varepsilon+2\sqrt{\gamma_p})$. As both left and right hand side functions in the inequality in \eqref{eq:temp12} are increasing functions in $\gamma_r$, we obtain that if
\begin{equation*}
\varepsilon(\varepsilon+2\sqrt{\gamma_p})\left(\sqrt{\gamma_p}+\varepsilon+3(\sqrt{\gamma_p}+\varepsilon+1)(\sqrt{\gamma_p}+\varepsilon+\varepsilon(\varepsilon+2\sqrt{\gamma_p})+\gamma_p)\right)^2< \gamma_p^2,  
\end{equation*}
then $\phi(1)-\phi(2)>0$ holds for all $\gamma_r \in (0, \varepsilon(\varepsilon+2\sqrt{\gamma_p}))$.
Finally, define, 
\begin{align*}
g(\varepsilon) := \varepsilon(\varepsilon+2\sqrt{\gamma_p})\left(\sqrt{\gamma_p}+\varepsilon+3(\sqrt{\gamma_p}+\varepsilon+1)(\sqrt{\gamma_p}+\varepsilon+\varepsilon(\varepsilon+2\sqrt{\gamma_p})+\gamma_p)\right)^2-\gamma_p^2.
\end{align*}
Then as $g(0)=-\gamma_p^2<0$ and $g$ is continuous in $\varepsilon$, there exists $\varepsilon_0>0$ such that $g(\varepsilon)<0$ for every $0<\varepsilon<\varepsilon_0$. To conclude, for every $\gamma_r\in ( 0,\varepsilon_0(\varepsilon_0+2\sqrt{\gamma_p}))$, we have both $\gamma_r+\gamma_p<\varepsilon_0(\varepsilon_0+2\sqrt{\gamma_p})+\gamma_p=(\varepsilon_0+\sqrt{\gamma_p})^2=\lambda_1^2$ and $\phi(1)>\phi(2)$. 

\textbf{Statement (iii):} For simplicity, we show the proof of statement (iii) again for $d=1$, $\sigma^2=1$. As in (ii), the general case is proven in an analogous way.  Fix $\lambda_1>0$. We need to show that there exists $\varepsilon \in (0, \lambda_1^2)$ such that for $\gamma_p=\lambda_1^2-\varepsilon$, we have $\phi(1)>\phi(2)$ for every $\gamma_r \in (0, \varepsilon)$. As in (ii) we have  
\begin{equation}\label{eq:temp2}
\gamma_r\left(\lambda_1+3(\lambda_1+1)(\lambda_1+\gamma_r+\gamma_p)\right)^2<\lambda_1^2(\gamma_r+\gamma_p)\Rightarrow \phi(1)-\phi(2)>0.  
\end{equation}
Substituting $\gamma_p=\lambda_1^2-\varepsilon>0$ in \eqref{eq:temp2}, we get 
$$
\gamma_r\left(\lambda_1+3(\lambda_1+1)(\lambda_1+\gamma_r+\lambda_1^2-\varepsilon)\right)^2<\lambda_1^2(\gamma_r+\lambda_1^2-\varepsilon)\Rightarrow \phi(1)-\phi(2)>0.  
$$
As in (ii) it is enough to show that 
$$
\varepsilon\left(1+3(\lambda_1+1)^2\right)^2<\lambda_1^2 - \varepsilon,
$$
for $\varepsilon>0$ small enough. That is, however, true, for every
\begin{align*}
    \varepsilon<\frac{\lambda_1^2}{\{ 1+3(\lambda_1+1)^2 \}^2 + 1},
\end{align*}
concluding the proof.
\end{proof}

\begin{proof}[Proof of Theorem \labeltheofour]
    We assume in the proof that $\gamma_p \in (0, 1]$, the case $\gamma_p \in (1, \infty)$ being proven analogously. As described in the main text, direct computation shows that $f_n(\tau_{j, n}) \rightarrow_p \lambda_j$ for all $j = 1, \ldots, d$. Similarly, one shows that
    \begin{align*}
    f_n(\tau_{j_n, n}) \rightarrow_p \frac{1}{2} \sigma^2 \{ F^{-1}_{\gamma_p}(q) - 1 - \gamma_p \},
    \end{align*}
    for all sequences $j_n > d$, such that $j_n/(p_n + r_n) \rightarrow 1 - q \in [0, 1]$. Consequently, plugging in to equation (2) in the main text, we obtain that $h_{j, n} \rightarrow_p \gamma_r \sigma^2$ for all $j = 1, \ldots, d$ and that 
    \begin{align}\label{eq:Aq_limit}
            h_{j_n, n} \rightarrow_p \frac{\gamma_r \sigma^2 A(q)(A(q) + 2\gamma_p + 2\gamma_r)}{2(\gamma_p + \gamma_r)(A(q) + 2)},
    \end{align}
    for all sequences $j_n > d$, such that $j_n/(p_n + r_n) \rightarrow 1 - q \in [0, 1]$ as $n \rightarrow \infty$, where $A(q) :=  F^{-1}_{\gamma_p}(q) - 1 - \gamma_p$. Hence, we have $h_{j + 1, n} - h_{j, n} \rightarrow_p 0$ for all $j = 1, \ldots, K$ such that $j \neq d$. For the noise part, this follows from the fact that two consecutive indices $j, j + 1$ correspond asymptotically to the same value of quantile index $q$. To compute the gap for the true signal dimension $d$, it is sufficient to plug in the value $A(1) = 2 \sqrt{\gamma_p}$ to the right-hand side of \eqref{eq:Aq_limit}.
\end{proof}

\begin{proof}[Proof of Lemma \labellemone]
The proof of Lemma~\labellemone \ is analogous to the proof of Theorem~\labeltheotwo \ and is thus omitted. 
\end{proof}

\bibliographystyle{chicago}
\bibliography{references}

\end{document}